\documentclass{article}
\usepackage[francais]{babel}
\usepackage{amsfonts}
\usepackage{xypic}

\def\id{{\rm id}}
\let\ra=\rightarrow
\let\Ra=\Rightarrow
\def\H{\it Hom}
\def\D#1{\Delta[#1]}
\def\cro#1{[#1]}

\let\ph=\varphi
\def\A#1{ \cro {#1}\times \cro n}
\def\DA#1{\D{#1}\times\D n}
\def\S{\cal S}
\def\Sr{\S_{\rm reg}}
\let\Bbb\mathbb
\begin{document}

\title{Ensembles simpliciaux r\'eguliers}
\author{Michel Zisman}
\date{10 septembre 2009}
\maketitle

{\bf R\'ESUM\'E}. On d\'efinit une sous-cat\'egorie int\'eressante $\Sr$ de la cat\'egorie des ensembles simpliciaux, dont les objets sont appel\'es {\it r\'eguliers}. Cette cat\'egorie, ainsi que sa sous-cat\'egorie 
${\cal S}_{f-{\rm reg}}$ dont les objets sont les ensembles simpliciaux r\'eguliers et finis, ont de bonnes propri\'et\'es de stabilit\'e par limite et union. La cat\'egorie ${\cal S}_{f-{\rm reg}}$ est cart\'esienne ferm\'ee, en contraste de celle des enembles simpliciaux finis qui n'est pas cart\'esienne ferm\'ee.

Le but de cette note est de r\'epondre \`a une question d'Andr\'e Joyal, qui m'a \'et\'e pos\'ee par Georges Maltsiniotis. D\'esignons par 
$\H$ le Hom interne de la cat\'egorie $\S$ des ensembles simpliciaux, et soient $U$ et $X$ deux ensembles simpliciaux finis (i.e. qui n'ont qu'un nombre fini de simplexes non d\'eg\'en\'er\'es) : dans ces conditions $\H (U,X)$ est-il fini ?
Curieusement, la r\'eponse est n\'egative en g\'en\'eral. Le premier contre-exemple, d\^u \`a  Jacob Lurie\footnote{au cours d'une conversation avec A. Joyal et G.Maltsiniotis}, consiste \`a prendre $U = \D1$ et  $X = \D3/\dot\D3$. Je n'ai pas r\'eussi jusqu'\`a pr\'esent \`a caract\'eriser les ensembles simpliciaux finis $X$ pour lesquels la r\'eponse est oui pour tous les ensembles simpliciaux finis $U$, mais la r\'eponse est oui pour une large famille d'ensembles simpliciaux, que l'on appellera ensembles simpliciaux {\it r\'eguliers} (confer 1.3.), en vertu du th\'eor\`eme suivant :

{\bf Th\'eor\`eme 1} : Si $X$ est un ensemble simplicial r\'egulier de dimension finie, alors,  pour tout ensemble simplicial fini $U$,
l'ensemble simplicial
$\H (U,X)$ est aussi de dimension finie.

(On dit qu'un ensemble simplicial est de {\it dimension finie} si le degr\'e de ses simplexes non d\'eg\'en\'er\'es est born\'e).

En particulier si $X$ est r\'egulier et fini, alors, pour tout ensemble simplicial fini $U$, l'ensemble simplicial $\H (U,X)$ est fini.

Nous verrons aussi que les sous-cat\'egories pleines $\Sr$ de $\S$ form\'ees par ces ensemble simpliciaux et ${\cal S}_{f-{\rm reg}}$ form\'ees par ceux qui sont en plus finis,
poss\`edent de bonnes propri\'et\'es de stabilit\'e. 

Je remercie Georges Maltsiniotis pour ses encouragements, le soin avec lequel il a relu de pr\'ec\'edents \'etats du manuscrit, et ses suggestions pour faciliter la lisibilit\'e du texte. Je remercie aussi  le referee dont les questions judicieuses m'ont permis de compl\'eter utilement  certains \'enonc\'es et de rendre plus agr\'eable la lecture de quelques d\'emonstrations.

{\bf 1. Pr\'eliminaires.}

{\bf 1.1.} Rappelons que les $p$-simplexes de $\H(U,X)$ sont les morphismes $f : \D p \times U \ra X$ de $\S$. Il en r\'esulte imm\'ediatement que l'on a 
$$\H ({\rm colim}_\alpha U_\alpha, X) 
               = {\rm lim}_\alpha\H( U_\alpha, X )  .$$
 Comme tout ensemble simplicial $U$ est colimite d'une famille de 
  $\D n_x$ (indic\'ee par les simplexes non d\'eg\'en\'er\'es $x$ de $U$) on voit que pour d\'emontrer que quel que soit $U$ fini, $\H (U,X)$  est fini pour un certain $X$, il suffit de le v\'erifier pour les $U = \D n$, et on se limitera \`a ce cas dans toute la suite. 
  
{\bf 1.2.} Voici quelques notations utiles dans la suite. Les objects de la cat\'egorie simpliciale $\Delta$ sont d\'esign\'es comme d'habitude par des entiers entre crochets. Les morphismes faces et d\'eg\'en\'erescences sont not\'es respectivement $\partial_i$ et 
$\sigma_i$.
Soient $X$ un ensemble simplicial et $\ph : \cro p\ra \cro q$ un morphisme de $\Delta$. On note  
$X(\ph) : X_q \ra X_p$ 
l'application associ\'ee par la structure simpliciale de $X$. Soit
 $\ph (i,i+r) :\cro 1\ra\cro p$ le morphisme d\'efini par $0\mapsto i$ et 
$1 \mapsto i+r$. 

{\bf 1.3.} Introduisons maintenant quelques d\'efinitions.  
  
On dit qu'un ensemble simplicial $X$ est 
 {\it fortement r\'egulier} si pour tout simplexe $x$ non d\'eg\'en\'er\'e de $X$, les faces $d_ix$ sont aussi non d\'eg\'en\'er\'ees. 
 
Une {\it ar\^ete \'el\'ementaire} d'un $n$-simplexe $ x$ d'un ensemble simplicial $X$ est un
1-simplexe $y$ de $X$ \'egal \`a $X(\ph(i,i+1))x$ pour un certain $i$, 
$0\leq i\leq n -1$.

On dit qu'un ensemble simplicial v\'erifie la propri\'et\'e 
$P_r$ si, \'etant donn\'e un simplexe $x$ de $X$ tel que  $X(\ph(i,i+r))x$ soit d\'eg\'en\'er\'e, alors il existe $y\in X$
tel que $x = s_{i+r-1}\ldots s_iy$.

On dit qu'un ensemble simplicial est {\it r\'egulier} s'il v\'erifie la propri\'et\'e $P_1$. Comme le laisse entendre la terminologie propos\'ee, nous verrons qu'un ensemble simplicial fortement r\'egulier est r\'egulier.

{\bf 1.4. Quelques propri\'et\'es \'el\'ementaires}

{\bf 1.4.1.}
On v\'erifie facilement le r\'esultat suivant :

(*) Soit  $x = s_{i+r-1}s_{i+r-2}\ldots s_iy $ un simplexe de $X$. Alors
l'ar\^ete  $X(\ph(i,i+r))x$ est  d\'eg\'en\'er\'ee.

On en d\'eduit imm\'ediatement le
 
{\bf Lemme 1} : Soit $y\in X$ un $q$-simplexe, soient $\alpha : \cro p \ra \cro q$ un morphisme surjectif, et $x = X(\alpha)y$. Alors 
$X(\ph(i,i+1))x$ est non d\'eg\'en\'er\'ee pour au plus $q$ valeurs de $i$.

En effet $X(\alpha) = s_{i_1}s_{i_2}\cdots s_{i_{p-q}}$ avec 
$i_1>i_2>\cdots>i_{p-q}$.

Lorsque $X$ est fortement r\'egulier, le r\'esultat $(*)$ poss\`ede une r\'eciproque :

{\bf Lemme 2} : Un ensemble simplicial fortement r\'egulier v\'erifie la propri\'et\'e $P_r$ pour tout $r>0$ (en particulier, il est r\'egulier). 
R\'eciproquement si un ensemble simplicial satisfait aux conditions 
$P_1$ et $P_2$, il est fortement r\'egulier.

{\it D\'emonstration} : Pour la premi\`ere partie, on proc\`ede par r\'ecurrence sur le degr\'e $p$ des simplexes, le cas $p=1$ \'etant tautologique. 
Soit $x\in X_p$ un $p$-simplexe tel que $X(\ph(i,i+r))x$ soit d\'eg\'en\'er\'e. Puisque $X$ est fortement r\'egulier, $x$ est donc d\'eg\'en\'er\'e, disons 
$x = s_ky$ pour un certain $k\leq p-1$. Il vient donc $X(\ph(i,i+r))x = X(\sigma_k\circ\ph(i,i+r))y$. Les relations de commutation 
$$   \sigma_k\circ\ph (i,i+r)   
=\left\{\matrix{
\ph(i-1,i+r-1) & {\rm si}\quad  k<i  \hfill    \cr
\ph(i,i+r-1)\hfill  & {\rm si}\quad i\leq k<i+r \cr
\ph(i,i+r) \hfill   & {\rm si}\quad k\geq i+r \hfill .   \cr
}\right.$$
permettent d'\'ecrire suivant les cas
 $$X(\ph(i,i+r))x =
\left\{\matrix{
X(\ph(i-1,i+r-1))y  \cr
X(\ph(i,i+r-1))y\hfill   \cr
X(\ph(i,i+r))y \hfill   .   \cr
}\right.$$

Puisque cette ar\^ete est d\'eg\'en\'er\'ee, l'hypoth\`ese de r\'ecurrence montre qu'il existe un $z$ de degr\'e $p-2$ tel que, selon les valeurs de $k$, on a 
$$y = 
\left\{\matrix{
s_{i+r-2}\cdots s_{i-1}z\cr
s_{i+r-2}\cdots s_iz\hfill \cr
s_{i+r-1}\cdots s_iz\hfill   .  \cr
}\right.$$
Mais alors il vient
$$x=s_ky =
\left\{\matrix{
s_{i+r-1}\cdots s_is_kz\hfill\cr
s_{i+r-1}\cdots s_iz\hfill \cr
s_{i+r-1}\cdots s_is_{k-r}z\hfill   .  \cr
}\right.$$

Supposons maintenant que $X$ v\'erifie $P_1$ et $P_2$. Soit $x\in X$ et supposons que $d_k x$ est d\'eg\'en\'er\'e, disons
$d_k x = s_l y$. Alors (confer $(*)$)
$X(\ph(l,l+1))d_k x$ est d\'eg\'en\'er\'e. Or on a 
$$X(\ph(l,l+1))d_k = X(\partial_k\circ\ph(l,l+1))$$
et on v\'erifie les \'egalit\'es suivantes :
$$\partial_k\circ\ph(l,l+1)
=\left\{\matrix{
\ph (l+1,l+2) & {\rm si} \quad k\leq l \hfill \cr
\ph (l,l+2) \hfill    &{\rm si} \quad k = l+1 \hfill  \cr
\ph (l,l+1) \hfill   &{\rm si}    \quad k>l+1 \hfill . \cr
}\right.$$

L'hypoth\`ese implique donc que, dans le premier cas, on a 
$x = s_{l+1}z $, dans le second on a $x = s_{l+1}s_l z$ et dans le troisi\`eme $x=s_l z$ pour un certain $z$. Dans tous les cas, $x$ est d\'eg\'en\'er\'e. 

{\bf Lemme 3} : Soit $X$ un ensemble simplicial.  Les trois propri\'et\'es suivantes sont \'equivalentes :
\begin{enumerate}
\item[(i)] $X$ est r\'egulier.
\item[(ii)] Les ar\^etes \'el\'ementaires d'un simplexe non d\'eg\'en\'er\'e de $X$ sont toutes non d\'eg\'en\'er\'ees.
\item[(iii)] Un simplexe $x$ de $X$ est non d\'eg\'en\'er\'e si et seulement si toutes ses ar\^etes \'el\'ementaires sont non d\'eg\'en\'er\'ees.
\end{enumerate}

{\it D\'emonstration} : (i) $\Ra$ (ii) par d\'efinition m\^eme de r\'egulier. (iii) n'est autre que (ii) \`a laquelle on ajoute la propri\'et\'e vraie sans restriction sur $X$ \`a savoir qu'un d\'eg\'en\'er\'e poss\`ede toujours une ar\^ete \'el\'ementaire d\'eg\'en\'er\'ee (confer $(*)$). Reste \`a montrer (ii)
$\Ra$ (i) et pour cela que si $X$ ne satisfait pas \`a (i), il existe un simplexe non d\'eg\'en\'er\'e de $X$ dont une ar\^ete \'el\'ementaire est
d\'eg\'en\'er\'ee. L'hypoth\`ese dit qu'il existe un $z\in X_p$ et un $i$ tel que 
$X(\ph(i,i+1))z$ est d\'eg\'en\'er\'e, mais qu'il n'existe aucun $y$ tel que 
$z = s_i y$. \'Ecrivons $ z = X(\phi )x$ avec $x$ non d\'eg\'en\'er\'e et $\phi$ surjective : $\phi(i)$ et $\phi(i+1)$ sont donc soit \'egaux, soit deux entiers successifs. Le premier cas ne peut se pr\'esenter car il impliquerait  l'existence d'un   $\psi$ tel que $\phi = \psi\circ\sigma_i$, et donc on aurait $z = s_iX(\psi)x$ en contradiction avec l'hypoth\`ese. Reste donc le second. Posons  $\phi (i) = j$. Comme on a 
$X(\ph(i,i+1))z = X(\ph(j,j+1))x$ 
et que $x$ est non d\'eg\'en\'er\'e, la d\'emonstration est achev\'ee. 

\medskip

La propri\'et\'e (ii) du lemme pr\'ec\'edent est tr\`es commode pour reconna\^\i tre un ensemble simplicial r\'egulier. Par exemple, si $n\geq 2$,  le quotient de $\D n$ par son ar\^ete $0n$, qui n'est \'evidemment pas fortement r\'egulier, est r\'egulier. De m\^eme on d\'emontre sans peine la proposition suivante :

{\bf Proposition 1} : La cat\'egorie $\Sr$ est stable par 
limites et sommes. Un sous-objet (dans $\S$) d'un objet de $\Sr$ est dans $\Sr$. Si $X_a$ est une famille de sous-ensembles simpliciaux de $X$, et si chaque $X_a$ est r\'egulier, la r\'eunion $\bigcup X_a$ l'est aussi. Un nerf est toujours r\'egulier.  

(Pour montrer par exemple la stabilit\'e de la cat\'egorie $\Sr$ par produits, il suffit d'utiliser la d\'efinition ; pour montrer qu'un sous ensemble simplicial d'un ensemble simplicial r\'egulier est r\'egulier, on utilise la propri\'et\'e (ii) et le fait qu'un simplexe d'un sous-ensemble simplicial est d\'eg\'en\'er\'e si et seulement si il est d\'eg\'en\'er\'e dans l'ensemble simplicial tout entier ; la stabilit\'e par limites r\'esulte de ces deux r\'esultats.)

\medskip
{\bf 2. Un crit\`ere de d\'eg\'en\'erescence}.

{\bf 2.1.} Rappelons que les $w$-simplexes de $\DA p$ sont les fonctions croissantes $\cro w \ra\A p$. Si la deuxi\`eme coordonn\'ee reste constante, on dira que le simplexe est {\it horizontal}, et si la premi\`ere coordonn\'ee reste constante, on dira qu'il est {\it vertical}.

Appelons {\it chemin} $a$ du r\'eseau $\A p$  un simplexe non d\'eg\'en\'er\'e, de longueur maximale pour origine et extr\'emit\'e fix\'ees. G\'eom\'etriquement, le chemin explicite les valeurs successives de 
$a : \cro w \ra \A p$, deux valeurs successives ayant toujours soit m\^eme abscisse  soit m\^eme ordonn\'ee. Notons $\cal C$ l'ensemble des
$\bigl({{n+p}\atop p}\bigr)$
chemins maximaux  i.e. ceux qui relient $(0,0)$ \`a $(p,n)$, et 
faisons la remarque, triviale mais utile, que si $c$ est un chemin maximal passant par le point de coordonn\'ees $(i,j)$, alors on a $c(i+j) = (i,j)$. 

Pour se donner 
un $p$-simplexe $f$ de $\H(\D n, X)$, donc un morphisme 
$f : \DA p \ra X$, il suffit de se donner les $(p+n)$-simplexes 
$z_a = f(a)$ de $X$ lorsque $a$ parcourt $\cal C$, ces donn\'ees \'etant soumises aux relations naturelles qui expriment que ``sur l'intersection $a\cap b$ de deux chemins maximaux, $z_a$ et $z_b$ co\"\i ncident"
\footnote{Confer par exemple  P. Gabriel and M. Zisman {\it Calculus of fractions and homotopy theory}, Ergebnisse der Mathematik, Band 35, Chapter II, 5.5.}.
On se propose dans ce paragraphe de donner un crit\`ere qui exprime que $f$ est en fait un d\'eg\'en\'er\'e d'un $(p-1)$-simplexe $g$. Commen\c cons par le diagramme commutatif suivant.

{\bf 2.2.} \'Etant donn\'e un chemin maximal $a$ de $\A p$ et un entier positif ou nul $k<p$, il existe un unique entier $t$ tel que l'on a
{\arraycolsep2pt
\begin{eqnarray*}
a(k+t)     &  =& (k,t)\hfill\\
a(k+t+1) &  =& (k+1,t)\hfill
\end{eqnarray*}}
Ayant ainsi fix\'e $t$, on d\'efinit un chemin maximal $m$ de $\A{p-1}$
en posant 
$$m(j) =
\left\{\matrix{
a(j)              & {\rm si} \quad j\leq k+t \hfill\cr
a(j+1)-(1,0) & {\rm si}  \quad j\geq k+t \hfill\cr
}\right.$$
et on v\'erifie sans peine que le diagramme
$$\xymatrix{
\cro{p+n}\ar[r]^a\ar[d]_{\sigma_{k+t}} & \A p\ar[d]^{\sigma_k\times\id}\\
\cro{p+n-1}\ar[r]_m &\A{p-1}
}$$

commute. Mais alors si $g : \DA {p-1}  \ra X$ est un $(p-1)$-simplexe de $\H (\D n, X)$ et si $h$ d\'esigne le $p$-simplexe d\'eg\'en\'er\'e 
$s_kg$, il vient : 
$$h(a) = s_{k+t}g(m).$$
D'apr\`es l'assertion $(*)$ de 1.4.1., l'image par $h$ de l'ar\^ete 
$((k,t), (k+1,t))$ est d\'eg\'en\'er\'ee. Introduisons, pour d\'esigner ce ph\'enom\`ene, les d\'efinitions suivantes :  

Soit $k$ un entier, $0\leq k\leq p-1$. On dit qu'un $p$-simplexe
$f$ de $\H(\D n, X)$ est {\it $k$-presque d\'eg\'en\'er\'e} si pour tout $0\leq j\leq n$, le 1-simplexe de $X$ \'egal
au compos\'e
$\xymatrix{
\cro 1\ar[r] & \A p\ar[r]^-f & X
}$
(o\`u la premi\`ere fl\`eche est d\'efini par $0 \mapsto (k,j)$,   
$1 \mapsto(k+1,j)$)
est d\'eg\'en\'er\'e. On dira qu'il est {\it presque-d\'eg\'en\'er\'e} s'il existe un $k$ tel qu'il est $k$-presque d\'eg\'en\'er\'e.

 Nous avons donc constat\'e qu'un $p$-simplexe d\'eg\'en\'er\'e de $\H (\D n, X)$ est presque d\'eg\'en\'er\'e. Nous verrons que moyennant des conditions sur $X$, cet \'enonc\'e poss\`ede une r\'eciproque. 

{\bf 2.2.1.}
Nous allons dans ce but pr\'eciser
la forme des chemins maximaux dans la tranche verticale limit\'ee par les points d'abscisse $k$ et $k+1$. Pour tout triplet $\alpha, t,  \beta$ d'entiers avec $0\leq \alpha\leq t\leq\beta\leq n$, 
soit
$b_{\alpha t\beta}$ l'unique chemin du r\'eseau $\A p$
passant par les points
$(k,\alpha), (k,t), (k+1,t)$ et $(k+1,\beta)$
et soit $l_{\alpha \beta}$ l'unique chemin du r\'eseau $\A {p-1}$ passant par les points $(k,\alpha)$ et $(k,\beta)$. Si $k = 0$,  on prend 
$\alpha = 0$, si $k=p -1$, on prend $\beta = n$.

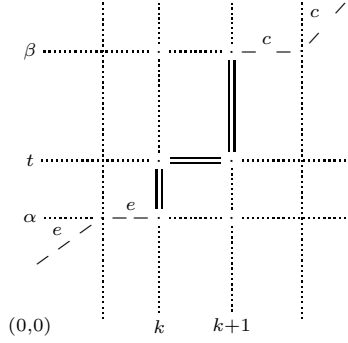
$$\def\objectstyle{\scriptstyle}
\xymatrix @-1pc{
&&&&&\\
\beta\ar@{.}[rr]&&.\ar@{.}[u]\ar@{.}[r]\ar@{.}[dd]&.\ar@{.}[u]\ar@{--}[r]^c&.\ar@{--}[ur]^c&\\
&&&&&&\\
t\ar@{.}[rr]&&.\ar@{=}[r]&.\ar@{=}[uu]\ar@{.}[rr]\ar@{.}[d]&&\\
\alpha\ar@{.}[r]&.\ar@{--}[r]^e&.\ar@{=}[u]\ar@{.}[r]&.\ar@{.}[rr]&&\\
\ar@{--}[ur]^e&&&&&\\
(0,0)&\ar@{.}[uuuuuu]&k\ar@{.}[uu]&k+1\ar@{.}[uu]&\ar@{.}[uuuuuu]&&
}$$

\centerline{\small en pointill\'e : le r\'eseau $\A p$, \  en trait double : 
  $b_{\alpha t \beta}$,  \  en tirets : $e$ et $c$.}

On introduit aussi les  ensembles de chemins $E_\alpha$ qui joignent
les points (0,0) \`a $(k,\alpha)$ et qui {\it arrivent horizontalement} en 
$(k,\alpha)$ (i.e. qui passent par $(k-1,\alpha)$ ;
ensemble vide si $k=0$), et $C_\beta$ qui joignent 
$(k+1,\beta)$ \`a $(p,n)$ et qui {\it partent horizontalement} de 
$(k+1,\beta)$  (i.e. qui passent par $(k+2,\beta)$ ;
ensemble vide si $k = p-1$). Enfin \`a $c\in C_\beta$, on associe le chemin $c^-$ de $\A{p-1}$, d\'efini par $c^-(i) = c(i)-(1,0)$.

Ces notations \'etant fix\'ees, il est clair que tout chemin maximal 
$a\in {\cal C} $ s'\'ecrit d'une et d'une seule mani\`ere comme une somme
$a = e+b_{\alpha t \beta} + c$, avec $e\in E_\alpha$ et $c\in C_\beta$, 
l'addition $+$ d\'esignant la concat\'enation : $a(i) = e(i)$ pour 
$i\leq k+\alpha$, $a(i) = b_{\alpha t \beta}(i-k -\alpha )$ pour
$k + \alpha \leq i \leq k+ 1 + \beta$ et $a(i) = c(i-k-1-\beta)$ pour
$k+1+\beta\leq i$.
Par ailleurs, la relation de 2.2. s'\'ecrit maintenant
$$h(e + b_{\alpha t \beta} + c) = s_{k + t} g(e + l_{\alpha\beta} + c^-)$$
pour $h = s_kg$.

{\bf 2.2.2. Le lemme principal.}

{\bf Lemme 4} : Soient $X$ un ensemble simplicial r\'egulier, et $f$ un simplexe $k$-presque d\'eg\'en\'er\'e de $\H(\D n, X)$. Alors il existe un simplexe $g$ tel que $f = s_kg$. Les simplexes presque d\'eg\'en\'er\'es de 
$\H(\D n, X)$ sont donc d\'eg\'en\'er\'es.

{\it D\'emonstration} : Soit $f$ un $p$-simplexe de $\H(\D n, X)$  presque d\'eg\'en\'er\'e.
Il existe donc un entier $k$, $0\leq k\leq p-1$,
 tel que pour tout $0\leq j\leq n$, le 1-simplexe de $X$ \'egal
au compos\'e
$\xymatrix{
\D 1\ar[r] & \D p\times \D n \ar[r]^-f & X
}$
(o\`u la premi\`ere fl\`eche est d\'efinie par $0 \mapsto (k,j)$,   
$1 \mapsto(k+1,j)$)
est d\'eg\'en\'er\'e. Soit $a$ un chemin maximal que l'on \'ecrit 
$a = e + b_{\alpha t \beta} + c$. L'ar\^ete $X(\ph(k + t,k+t +1))f(a)$ est donc d\'eg\'en\'er\'ee. Puisque $X$ est r\'egulier, cela signifie que l'on a
$$f(a) = s_{k + t}y_t$$
pour un certain $y_t \in X$. 
Si $\alpha<t \leq\beta$,
choisissons maintenant $a' =  e + b_{\alpha (t-1) \beta} + c$. Il vient de m\^eme $f(a') =s_{k + t -1} y_{t-1}$. 
Comme $a$ et $a'$ ne diff\`erent que sur le carr\'e form\'e par les deux verticales du r\'eseau, d'abscisse $k$ et $k+1$ d'une part et les deux horizontales d'ordonn\'ee $t-1$ et $t$ d'autre part, les relations de compatibilit\'e indiquent que l'on a 
$d_{k + t}f(a) = d_{k + 1 + t -1}f(a')$
ce qui impose 
$$y_ t  = y_{t-1}.$$
Les $y_t$ sont donc ind\'ependants de $t$, et tous \'egaux \`a 
$d_{k + \beta + 1}f(e + b_{\alpha\beta\beta} + c)$, ne d\'ependant que de 
$e, c, \alpha, \beta$.
On les notera $g( e + l_{\alpha\beta} + c^-)$. Remarquons que si on a $\alpha = \beta$, alors le chemin $b_{\alpha\beta\beta}$  est r\'eduit \`a une ar\^ete
et le chemin $l_{\alpha\beta}$ \`a un point.

Ainsi, \`a tout chemin maximal $e + l_{\alpha\beta} + c^-$ du 
r\'eseau $\A {p-1}$, nous avons associ\'e le $(p + n-1)$-simplexe
$g( e + l_{\alpha\beta} + c^-)$ de $X$.
Par construction, on a 
$f(e + b_{\alpha t \beta} + c) = s_{k + t} g(e + l_{\alpha\beta} + c^-)$
et il reste \`a v\'erifier que les relations de compatibilit\'e sont satisfaites par les  $g( e + l_{\alpha\beta} + c^-) $. Soient $m$ et $m'$ deux chemins maximaux du r\'eseau $\A {p-1}$ qui ne diff\`erent que sur le carr\'e passant par les deux points $(i,j)$ et $(i + 1, j + 1)$. On doit avoir  $d_{i + j +1}g(m) = d_{i + j + 1}g(m')$. Il y a quatre cas, selon que  
$i<k-1, i = k-1, i = k$ et $i> k$. Les deux cas extr\^emes sont \'evidents. Traitons par exemple le cas $i = k$. 
On a $m = e + l_{\alpha\beta} + c^-$ et 
$m' = e + l_{\alpha(\beta -1)} + c'^-$,
avec $c'(0) = (k+1, \beta -1), c'(1) = (k+2,\beta -1)$ et 
$c'(r + 1) = c(r)$ pour $r>1$.
Nous devons v\'erifier l'\'egalit\'e
$$d_{k + \beta }d_{k + \beta + 1}f(e + b_{\alpha\beta\beta} + c) =
d_{k + \beta }d_{k + \beta }f(e + b_{\alpha(\beta -1(\beta-1)} + c')$$
ce qui ne pose aucune difficult\'e : c'est exactement la compatibilit\'e des deux $( p+ n)$-simplexes qui figurent de part et d'autre de l'\'egalit\'e.

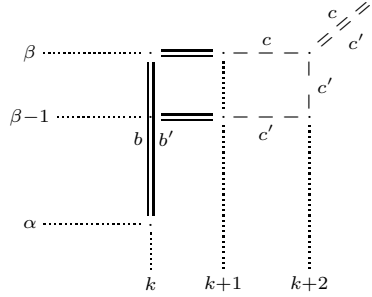
$$\def\objectstyle{\scriptstyle}
\xymatrix@-1pc{
&&&&&\\
 \beta\ar@{.}[rr]&&.\ar@{=}[r]&.\ar@{--}[r]^c&.\ar@2{--}[ur]^c_{c'}\\
 \beta -1\ar@{.}[rr]&&.\ar@{=}[r]&.\ar@{.}[u]  \ar@{--}[r]_{c'}&.\ar@{--}[u]_{c'}\\
 &&&&&\\
 \alpha\ar@{.}[rr]&&.\ar@{=}[uuu]^b_{b'}\\
 &&k\ar@{.}[u]&k+1\ar@{.}[uuu]&k+2\ar@{.}[uuu]
}$$

\centerline{\small on a \'ecrit $b$ au lieu de $b_{\alpha\beta\beta}$
et $b'$ au lieu de $b_{\alpha(\beta-1)(\beta-1)}$.}

\medskip
{\bf 3. Les th\'eor\`emes}.

{\bf 3.1.} La proposition suivante est un cas particulier du th\'eor\`eme 1 :

{\bf Proposition 2} : Si $X$ \ est un ensemble simplicial r\'egulier de dimension $q$, alors 
$\H (\D n,X)$ est de dimension au plus $(n+1)q$.

{\it D\'emonstration} : Soit $f$ un $p$-simplexe de $\H(\D n, X)$, i.e. un morphisme $\DA p \ra X$. Pour tout
$0\leq j\leq n$ , soit
 $x_j : \cro p\ra \cro p \times \cro n$ le $p$-simplexe de $\DA p$ d\'efini par $i \mapsto (i,j)$. Son image par $f$ est un $p$-simplexe de $X$. On dira que l'ar\^ete $((i,j),(i+1,j))$  du r\'eseau $\A p$ est {\it efficace} si l'ar\^ete \'el\'ementaire
$X(\ph(i,i+1))f(x_j)$ est non d\'eg\'en\'er\'ee.
Supposons maintenant que $X$ est de dimension finie $q$, et prenons $p>q$ : sous ces conditions  le simplexe 
$f(x_j)$ est d\'eg\'en\'er\'e ; le lemme 1 nous dit que ce simplexe poss\`ede au plus $q$ ar\^etes efficaces. Comme $j$ prend $n+1$ valeurs distinctes, le r\'eseau poss\`ede au plus $(n+1)q$ ar\^etes efficaces distinctes. Choisissons donc $p > (n+1)q$, ce qui nous assure de l'existence d'un entier $k$ tel que, pour tout $j = 0,\ldots n$, l'ar\^ete $((k,j),(k+1,j))$ est non efficace. En d'autres termes, le simplexe $f$ est presque d\'eg\'en\'er\'e. D'apr\`es le lemme 4, il est donc d\'eg\'en\'er\'e.

{\bf 3.1.1.} La borne $(n+1)q$ est la meilleure possible puisque 
$Hom(\D n,\D q)$ est de dimension $(n+1)q$.
 En effet les $p$-simplexes de $Hom(\D n,\D q)$ sont les applications croissantes $f :  \A p \ra \cro q$ et ce simplexe est d\'eg\'en\'er\'e si et seulement si $f$ prend les m\^emes valeurs sur deux colonnes voisines du r\'eseau $\A p$. Soit alors $f :  \cro {(n+1)q} \times \cro n \ra  \cro q$ 
 l'application d\'efinie par les \'egalit\'es suivantes, o\`u $j$ parcourt les valeurs $0\leq j\leq n$, et o\`u l'on a \'ecrit  $i = kq + a $ pour $i>0$, 
 $k\in\{0,\ldots,n\}$ et $1\leq a \leq q$ : 
on pose $f(0,j) = 0$  et, pour $i>0$,
 $$f(i,j) =
\left\{\matrix{
0 & {\rm si} \quad j<n-k \hfill\cr
a & {\rm si}  \quad j=n-k \hfill\cr
q & {\rm si}  \quad j>n-k . \hfill\cr
}\right.$$
On d\'efinit ainsi un $(n+1)q$-simplexe non d\'eg\'en\'er\'e de 
$Hom(\D n,\D q)$ qui est donc de dimension au moins $q$. Montrons que la dimension est exactement $(n+1)q$. Soit $f$ un simplexe non d\'eg\'en\'er\'e de degr\'e $p$. Pour tout 
$i \in \cro p$, il existe au moins un $j\in \cro n$ tel que 
$f(i,j)<f(i+1,j)$. Donc, si $\sigma(i)$ d\'esigne la somme
$f(i,0)+f(i,1)+\cdots+f(i,n)$ des valeurs prises par la fonction sur la $i$-\`eme colonne du r\'eseau, il vient 
$\sigma(i)<\sigma(i+1)$ et finalement 
$\sigma (p)\geq p$. Mais par ailleurs il est clair que l'on a, pour  
tout $i$, $\sigma(i) \leq (n+1)q$, puisqu'il y a $n+1$ termes dans chaque colonne du r\'eseau. Donc il vient 
$p\leq\sigma(p)\leq (n+1)q$.

{\bf 3.1.2.}{\it D\'emonstration du th\'eor\`eme 1}. On \'ecrit $U = {\rm colim }_u{\D {|u|}}_x$ 
o\`u $u$ parcourt l'ensemble des simplexes non d\'eg\'en\'er\'es de $U$,
et o\`u $|u|$ d\'esigne le degr\'e de $u$. L'\'egalit\'e donn\'ee dans (1.1) et la proposition 2, ainsi que les propri\'et\'es \'el\'ementaires de la dimension permettent de conclure. Plus pr\'ecis\'ement, soient $A$ et $B$ deux ensembles simpliciaux  de dimension finie. On a 
${\rm dim} \ ( A\times  B) = {\rm dim}\ A + {\rm dim }\ B$, et
si $A\subset B$ alors il vient  dim $A <$ dim $B$. Comme une limite finie n'est autre qu'un sous objet d'un produit fini, nous obtenons :

{\bf Th\'eor\`eme 1bis} : Soient $X$ un ensemble simplicial de dimension finie et $U$ un ensemble simplicial fini. On a  : 
$$  {\rm dim} \ Hom (U,  X)  \leq\sum_u (|u|+1).{\rm dim} \ X$$
o\`u, dans la somme, $u$ parcourt l'ensemble des simplexes non d\'eg\'en\'er\'es de $U$.

{\bf 3.2.} Partant d'un ensemble simplicial r\'egulier $X$ qu'en est-il de 
$\H(U, X)$ ? Le th\'eor\`eme 2 r\'epond \`a la question :

{\bf Th\'eor\`eme 2} : Soient $X$ un ensemble simplicial r\'egulier et $U$ un ensemble simplicial quelconque. Alors $\H(U,X)$ est r\'egulier.

{\it D\'emonstration} : Utilisant la remarque 1.1. et la proposition 1, nous voyons qu'il suffit de d\'emontrer le th\'eor\`eme dans le cas o\`u $U = \D n$.
Posons $Y = \H(\D n, X)$. Soit $f\in Y_p$, et supposons que  $Y(\ph(k,k+1))f$ est d\'eg\'en\'er\'e. Cela signifie qu'il existe un $h$ qui fait commuter le diagramme suivant :
$$\xymatrix{
\D 1 \times\D n\ar [rr]^{\ph(k,k+1)\times \id}\ar[dr]_{\sigma_0} && \D p \times \D n \ar [r]^-f & X\\
 & \D n \ar [urr]_h& 
}$$

Mais alors pour tout $j \in \{0,\cdots, n\}$,
les 1-simplexes $((k,j),(k+1,j))$ ne sont pas efficaces. Comme dans la d\'emonstration du th\'eor\`eme 1, cela signifie que $f$ est $k$-presque d\'eg\'en\'er\'e. Le lemme 4 implique
 qu'il existe $g$ tel que $f = s_kg$. Donc $Y$ v\'erifie  $P_1$.

Lorsque l'on se restreint aux ensembles simpliciaux finis, on obtient (confer proposition 1): 

{\bf Th\'eor\`eme 3} : La cat\'egorie ${\cal S}_{f-{\rm reg}}$ 
est stable par limites finies et cart\'esienne ferm\'ee. Elle est
aussi stable par sous-objets et sommes finies (et m\^eme r\'eunions
finies).

{\bf 4. Quelques compl\'ements}. Disons qu'un ensemble simplicial $X$ est fortement fini si $\H(U,X)$ est fini pour tout ensemble simplicial fini $U$.

On peut montrer que tous les quotients de $\D 2$ sont fortement finis.
Il est probable qu'on peut en d\'eduire que les ensembles simpliciaux finis de dimension 2 sont fortement finis. 

Notons $F_i$ la $i$-\`eme face de $\D q$.Voici un r\'esultat qui g\'en\'eralise l\'eg\`erement celui de Lurie annonc\'e au d\'ebut :

{\bf Proposition 3 }Soit $F \subset \dot\D q$  ($q\geq 3$) une r\'eunion de faces qui contient $F_a\cup F_{a+1}$ pour un certain $a$, $0<a<q-1)$. Alors $X = \D q/F$ n'est pas fortement fini.

{\it D\'emonstration} (Lurie) : Introduisons la fonction coupe d\'efinie sur 
$\Bbb Z$ et \`a valeurs dans $\Bbb N$
par 
$${\rm coupe}(i) = 
\left\{\matrix{
 0&{\rm si} & i\leq 0\cr
 i & {\rm si} & 0\leq i \leq q\cr
 q& {\rm si} & q\leq i \quad . \cr
 }\right.$$     
Soit $p>q$ un entier, soit $0\leq u\leq p$ et 
soit $z_u$ le $(p+1)$-simplexe de $\D q$, i.e. le morphisme $\cro {p+1}\ra \cro q$ de la cat\'egorie $\Delta $, d\'efini par 
$z_u(i) = {\rm coupe}(i-u+a)$. 
Il est clair  que $z_u\circ\partial_u$ ne prend pas la valeur $a$ et que $z_u\circ\partial_{u+1}$ ne prend pas la valeur 
$a+1$. Ainsi $d_uz_u$ et $d_{u+1}z_u$ sont des simplexes de $F$.
Les $p+1$ simplexes $\bar z_u$, images des pr\'ec\'edents dans le quotient $X$ v\'erifient donc les relations de compatibilit\'e 
$d_{u+1}\bar z_u = d_{u+1}\bar z_{u+1}$  et d\'efinissent ainsi un $p$-simplexe de $Hom(\D 1, X)$. Ce simplexe n'est jamais d\'eg\'en\'er\'e.
En effet (d'apr\`es 2.2., ou un raisonnement direct), la condition pour qu'un $p$-simplexe $f$ de $\H (\D1,Z)$, donn\'e par  $p+1$ simplexes 
$f_u$, $u \in \{0,\ldots,p\}$, de degr\'e $p+1$ de $Z$, soit le d\'eg\'en\'er\'e
$s_k g$  d'un simplexe $g$ donn\'e par des $g_v\in Z_p $, est  
$$ f_u = s_{k+1}g_u \quad {\rm si }\quad  u\leq k \quad {\rm et }\quad 
           f_u = s_kg_{u-1}\quad {\rm si }\quad  u> k .$$
Il suffit de montrer que  $\bar z_k$ n'est pas dans l'image de $s_{k+1}$ pour v\'erifier l'assertion. 
Or on a $z_u(u+1) = a+1$ et $z_u(u+2) =a+2$.
Si on avait $\bar z_k = s_{k+1}\bar y$  pour un certain $y$, alors
$z_k= y\circ \sigma_{k+1}$ prendrait la m\^eme valeur pour $k+1$ et 
$k+2$, ce qui est impossible. Mais alors $Hom(\D 1, X)$ est de dimension infinie.

{\bf Remarque}: Qu'en est-il de $\D q/F_i$ ? La question reste ouverte. 
Je ne sais pas non plus si l'on peut remplacer dans le th\'eor\`eme 2, r\'egulier par fortement r\'egulier. C'est vrai pour $X = \D q$, d'apr\`es un raisonnement analogue \`a celui de (3.1.1.) utilisant la croissance de 
$\sigma$.

\bigskip

\vtop to0pt{\small\noindent
Michel Zisman\\
Universit\'e Paris 7\\
 zisman@math.jussieu.fr\vss}

\end{document}